\documentclass[letterpaper, 10 pt, conference]{ieeeconf}  

\IEEEoverridecommandlockouts                              
\title{\LARGE \bf
Worst-Case Risk Quantification under Distributional Ambiguity\\using Kernel Mean Embedding in Moment Problem
}

\author{Jia-Jie Zhu, Wittawat Jitkrittum, Moritz Diehl, Bernhard Sch\"olkopf
\thanks{This project has received funding from the European Union’s Horizon 2020 research and innovation programme under the Marie Skłodowska-Curie grant agreement No 798321, the German Federal Ministry for Economic Affairs and Energy (BMWi) via eco4wind (0324125B) and DyConPV (0324166B), and by DFG via Research Unit FOR 2401. \jz{any changes to funding sources?}}
\thanks{Jia-Jie Zhu, Wittawat Jitkrittum, and Bernhard Sch\"olkopf are with the Empirical Inference Department, Max Planck Institute for Intelligent Systems, T\"ubingen, Germany.
        {\tt\small \{jzhu,wittawat, bs\}@tuebingen.mpg.de}}%
\thanks{Moritz Diehl is with the Department of Microsystems Engineering, University of Freiburg, Freiburg, Germany.
        {\tt\small moritz.diehl@imtek.uni-freiburg.de}}%
}

\usepackage[unicode=true]{hyperref}
\usepackage{graphicx}

\usepackage{amsthm}
\usepackage{amsmath}
\usepackage{amsfonts}
\usepackage{dsfont} 

\newcommand{\rkhs}{{\mathcal H}}
\newcommand{\kme}[2]{\sum_{i=1}^N{#1} \phi({#2})}
\newcommand{\seqin}[1]{\{{#1}\}_{i=1}^N} 
\newcommand{\indfn}{\mathds{1}}

\DeclareMathOperator*{\mxz}{\mathrm{maximize}}
\DeclareMathOperator*{\sjt}{\text{subject to }}

\newtheorem{lemma}{Lemma}

\newtheorem{prop}{Proposition}

\newtheorem*{remark}{Remark}

\usepackage[normalem]{ulem} 

\usepackage{marginnote}
\usepackage{xcolor}

\usepackage{subcaption}

\usepackage{comment}

\excludecomment{todo}

\newif\ifcomment
\commentfalse

\ifcomment
\newcommand{\wjsay}[1]{[\textcolor{orange!60!black}{\textbf{WJ:}} \textcolor{orange!60!black}{\textbf{#1}}]}
\newcommand{\jz}[1]{[\textcolor{blue!60!black}{\textbf{JZ:}} \textcolor{blue!60!black}{\textbf{#1}}]}
\else
\newcommand{\wjsay}[1]{}
\newcommand{\jz}[1]{}
\fi
\begin{document}

\maketitle
\thispagestyle{empty}
\pagestyle{empty}

\begin{abstract}
In order to anticipate rare and impactful events, we propose to quantify the worst-case risk under distributional ambiguity using a recent development in kernel methods --- the kernel mean embedding. Specifically, we formulate the generalized moment problem whose ambiguity set (i.e., the moment constraint) is described by constraints in the associated reproducing kernel Hilbert space in a nonparametric manner. We then present the tractable approximation and its theoretical justification. As a concrete application, we numerically test the proposed method in characterizing the worst-case constraint violation probability in the context of a constrained stochastic control system.  
%
%
\end{abstract}

\section{Introduction}\label{sec:intro}
\begin{figure}[t!]
        \centering
        \includegraphics[width=0.5\columnwidth]{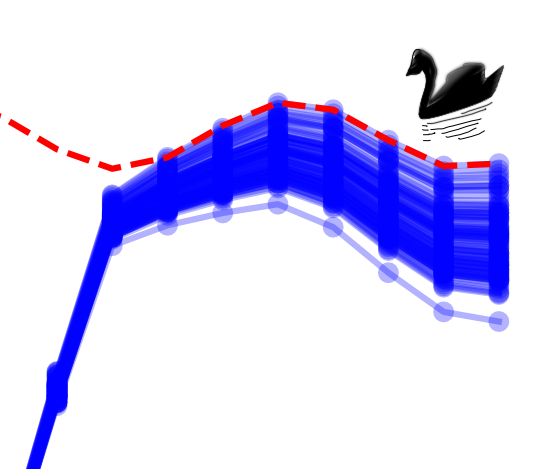}
        \caption{Rare constraint violation event as the black swan}
        \label{fig:swan}
\vspace{-0.5cm}
\end{figure}

We begin our discussion with the black swan metaphor, which was visited in a recent popular book \cite{talebBlackSwanImpact2007}.
It is argued that extremely rare events, like a black swan, may potentially have a huge impact on the underlying system. However, given limited historical data, statistical inference is prone to failures in predicting black swans. This issue is also relevant to optimization and control. 
Consider the illustrative Fig.~\ref{fig:swan} of a stochastic control system.
All the sampled state trajectories (denoted by $x_t$, in blue) satisfy some underlying state constraint, which is staying below the dashed curve (denoted by $c(x_t)\leq0$, in red).
If we are interested in estimating the constraint violation probability at a certain time from this empirical dataset, a naive Monte Carlo estimate yields 
$$
P{(c(X_t)>0)} \approx \sum_{i=1}^N\frac1N\indfn\{c(x_t^{(i)})>0\} = 0,
$$
i.e., zero violation probability. 
 This is sometimes referred to as the zero-count problem, also known as the silent evidence or inductivist turkey. However, just like the black swan metaphor, the constraint violation event may be rare but potentially impactful to the system. 

The key here is to capture the concept of distributional \emph{ambiguity}, i.e., we are uncertain about our knowledge of the probability measure itself.
For example, the classical Cantelli's inequality for a zero mean unit variance random variable $X\sim P$ states
$
P(X > c)\leq \frac{1}{1+c^2},
$
regardless of what distribution \(P\) is.
It is commonly argued in the literature of optimization and control, e.g., \cite{popescuSemidefiniteProgrammingApproach2005,farinaStochasticLinearModel2016}, that such bounds are pessimistic.
However, if $P$ is subject to perturbations, then the statistics such as mean and variance may become uncertain.
We refer to such uncertainty in distributions as \emph{ambiguity}.
For example, in machine learning applications such as \cite{sinhaCertifyingDistributionalRobustness2017a}, individual data samples themselves may be perturbed by an adversary or nature. 
Consequently, the original already-conservative bound may fail to hold.
From this point of view,
we must ask, ``how can nature perturb the empirical data distribution to hurt us the most?"

\subsubsection*{Generalized moment problem}
In the mathematics literature, this is the study of the moment problem such as the classical moment problems of Stieltjes, Hausdorff, and Hamburger a century ago. 
In a nutshell, a generalized moment problem can be abstractly written as
\begin{equation}
        \underset{P\in\mathcal C}{\mxz}\ \mathbb E _ P l(x),
\label{eq:moment}
\end{equation}
where \(\mathcal C\) is a set of probability measures. In this paper, we refer to \(\mathcal{C}\) as the \emph{ambiguity set}. 
For example, \(\mathcal C =\{P\mid \int \phi_i(x) \ dP(x) = m_i, i=1\dots p\}\) where \(\phi_i\)'s are some given moment functions and \(m_i\)'s are moments. 
This optimization problem searches for the \emph{worst-case risk} --- expected value of $l(x)$ under the \emph{worst-case distribution}.

The literature related to the moment problem is vast spanning over decades.
We now discuss a thread of modern computational approaches to deriving a bound (e.g., right-hand-side of Cantelli's inequality).
The early work of \cite{isiiSharpnessTchebychefftypeInequalities1962} contains the duality approach to the moment problem.
Subsequently, the author of \cite{shapiroDualityTheoryConic2001} proved strong duality using conic linear programming and the finite dimension reduction results from \cite{rogosinskiMomentsNonNegativeMass}.
Modern computational approaches such as \cite{lasserreBoundsMeasuresSatisfying2002,bertsimasMomentProblemsSemidefinite2000,bertsimasOptimalInequalitiesProbability2005} proposed SDP formulations for the moment problem under the assumption that the moment functions are polynomial-representable. 
Built on that work, a subsequent paper \cite{popescuSemidefiniteProgrammingApproach2005} used Choquet theory to study a few classes of non-degenerative distributions (e.g. unimodal) to reduce the conservativeness of the bound.
The authors of \cite{vandenbergheGeneralizedChebyshevBounds2007} showed that a worst-case probability problem admits exact SDP reformulation under the conditions that constraint set is a quadratic inequality and ambiguous distributions share the first two moments. 
Under similar ambiguity assumptions, the work of \cite{vanparysGeneralizedGaussInequalities2016} also used Choquet theory for the non-degenerative classes of distributions to derive exact reformulations under polyhedral constraints. 
Another related thread of work is the distributionally robust optimization (DRO) (e.g., \cite{delageDistributionallyRobustOptimization2010,erdoganAmbiguousChanceConstrained2006,mohajerinesfahaniDatadrivenDistributionallyRobust2018b}), as well as related approaches in data-driven robust optimization such as \cite{ben-talRobustSolutionsOptimization2013,bertsimasDatadrivenRobustOptimization2017}, which studies decision-making under distributional ambiguity. In a nutshell, it aims to find a decision after the worst-case distribution in \eqref{eq:moment} is assumed. 
Therefore, while this paper does not study DRO, the moment problem naturally becomes the inner optimization problem for DRO. Note the generalized moment problem is also referred to as the uncertainty quantification problem~\cite{owhadi2013optimal}.

At the core of those works is the description for the ambiguity set, e.g., the set of distributions $\mathcal C$ in \eqref{eq:moment}. 
We summarize a few example approaches taken in Table~\ref{table:ambiguity}. 
\begin{table}[!t]
        \caption{Examples of distributional ambiguity descriptions}
        \label{table:ambiguity}
        \centering
        \begin{tabular}{|c||c|}
        \hline
        \cite{bertsimasOptimalInequalitiesProbability2005,popescuSemidefiniteProgrammingApproach2005,vandenbergheGeneralizedChebyshevBounds2007,vanparysGeneralizedGaussInequalities2016}& Known first $p$ moments \\
        \hline
        \cite{erdoganAmbiguousChanceConstrained2006} & Prokhorov metric\\
        \hline
        \cite{ben-talRobustSolutionsOptimization2013} & $\phi$-divergence\\
        \hline
        \cite{mohajerinesfahaniDatadrivenDistributionallyRobust2018b} & Wasserstein metric\\
        \hline
        This paper & RKHS metric\\
        \hline
        \end{tabular}
\end{table}
To characterize ambiguity,
this paper summons a recent development in kernel methods, the kernel mean embedding. This framework can represent probability distributions as functionals in the reproducing kernel Hilbert space.
This elegant structural mapping allows us to cast worst-case probability measures as decision variables of optimization problems, as we shall introduce shortly.

The \emph{main contributions} of this paper are in using the kernel mean embedding for distributional ambiguity description, as well as the tractable approximate formulations (e.g. \eqref{eq:kme_mp},\eqref{eq:kmemp_prac}) and their theoretical justification (e.g. we prove a novel convergence result in Proposition~\ref{thm:cvgs}).
We then demonstrate how to use them in practice in the context of quantifying worst-case constraint violation probability of a stochastic control system.

\subsubsection*{Organization}
The paper is organized as follows.
Sec.~\ref{sec:kme} presents the background on the kernel mean embedding framework.
We then propose our main method in Sec.~\ref{sec:method}, whose proofs are deferred to 
later in Sec.~\ref{sec:theory}.
In Sec.~\ref{sec:experiment}, we analyze the proposed method in two numerical experiments: a synthetic worst-case probability problem and a stochastic control problem. 
The paper is concluded in Sec.~\ref{sec:conclusion}.

\subsubsection*{Notation}\label{sec:notation}
In this paper, $\mathcal X \subseteq \mathbb{R}^n$ denotes the domain or input space of interest, e.g., the state space of a system.
$\mathcal M$ denotes the vector space of all signed measures defined on $\mathcal B(\mathcal X)$, where $B(\mathcal X)$ is the Borel algebra.
$\mathcal P$ denotes the set of all probability measures on $\mathcal B(\mathcal X)$.
$\rkhs$ denotes the reproducing kernel Hilbert space associated with the kernel in the context.
$\int$ denotes the Stieltjes integral for both discrete and continuous distributions.
$M$ often denotes the number of data samples we have.
$\delta_{x}$ is the Dirac measure at $x$.
By embedding, we refer to the kernel mean embedding which we shall introduce next.

\hypertarget{kme}{%
\section{Preliminaries on Kernel mean embedding}\label{sec:kme}}
In this section, we give a brief overview of the kernel mean embedding --- a modern development in kernel methods --- which enables us
to describe distributional ambiguity non-parametrically. 
We start with the definition of a kernel. 
A symmetric function $k\colon \mathcal{X}
\times \mathcal{X} \to \mathbb{R}$ is said to be a positive definite kernel
if for any $N \in \mathbb{Z}^+$, any $\{ x_1, \ldots, x_N \} \subset
\mathcal{X}$, and any $\alpha_1,\ldots, \alpha_N \in \mathbb{R}$, we have
$\sum_{i=1}^N \sum_{j=1}^N \alpha_i \alpha_j k(x_i, x_j) \ge 0.$ We will
subsequently refer to a positive definite kernel  simply as a kernel.
Given a kernel $k$, the Moore-Aronszajn theorem \cite{berlinetReproducingKernelHilbert2011} guarantees
the existence of a Hilbert space $\rkhs$ of real-valued functions, and a map
$\phi\colon \mathcal{X} \to \rkhs$, such that $\langle f, \phi(x)
\rangle_\rkhs = f(x)$ for all $f \in \rkhs, x \in \mathcal{X}$, where
$\langle \cdot , \cdot \rangle := \langle \cdot , \cdot \rangle_\rkhs$
denotes the inner product on $\rkhs$. This property is known as the
reproducing property, and $\rkhs$ is known as the reproducing kernel Hilbert
space (RKHS) associated with the kernel $k$. 
The map $\phi$ is commonly known as a feature map. The reproducing property implies
that $k(x,x') = \langle \phi(x), \phi(x') \rangle$ for any $x, x' \in
\mathcal{X}$. 
Therefore, we interchangeably write $k(x, \cdot)$ and $\phi(x)$.

The kernel mean embedding \cite{smolaHilbertSpaceEmbedding2007} is a
technique that gives nonparametric representations of distributions in an RKHS. 
More precisely, given a probability distribution $P$ and a kernel $k$, the
(kernel mean) embedding of $P$ is defined as
\begin{equation}
\mu_P(\cdot) = \int  k(x, \cdot)  \ dP(x),
\label{eq:defkme}
\end{equation}
which is simply the expectation of the feature map under distribution $P$.
If $\mathbb{E}_{x \sim P} k(x,x) < \infty$, then $\mu_P$ exists and is in
$\rkhs$ \cite{smolaHilbertSpaceEmbedding2007}. One may also view $\mu_P \in \rkhs$
as a vectorial representation of $P$ in $\rkhs$. For example, if $\mathcal{X}=\mathbb{R}, \phi(x) :=
(x,x^2)^\top$, then $k(x,y) = xy + x^2 y^2$, and $\mu_P = \left(
\mathbb{E}_{x \sim P}[x], \mathbb{E}_{x \sim P}[x^2] \right)^\top $ which contains the first two moments of $P$.
More generally, if the kernel $k$ is \emph{characteristic}, then the kernel
mean map $P \mapsto \mu_P$ is injective so that $\mu_P$ is unique for any
distribution $P$ \cite{sriperumbudurUniversalityCharacteristicKernels2011}. An example of characteristic kernels
is the Gaussian kernel $k(x,y) = \exp\left( -\frac{ \| x-y\|^2_2}{2 \sigma^2}
\right)$ for any bandwidth $\sigma \in (0, \infty)$. This kernel
defines an infinite-dimensional feature map $\phi$.

The embedding can be used to define a distance between two
distributions,
Specifically, given two Borel probability measures $P,Q$, and a
characteristic kernel, one can show that $\| \mu_P - \mu_Q \|_\rkhs$ defines
a proper distance in the associated RKHS~$\rkhs$. 
This distance is also known as the maximum mean discrepancy (MMD) \cite{grettonKernelTwosampleTest2012}.
Note that $\| f \|_\rkhs := \sqrt{ \langle f, f\rangle_\rkhs} $ for any $f \in
\rkhs$. 
It follows that
\begin{multline*}
\| \mu_P - \mu_Q \|^2_\rkhs 
= \mathbb{E}_{x,x'\sim P} k(x,x') \\
- 2\mathbb{E}_{x\sim P} \mathbb{E}_{y \sim Q} k(x,y)
+ \mathbb{E}_{y,y' \sim Q} k(y, y'),
\end{multline*}
which can be estimated with the plug-in estimator using the samples from $P$
and $Q$ by virtue of the kernel trick \cite{grettonKernelTwosampleTest2012}.

This paper's idea is to use the embedding and the distance in RKHS to characterize the distributional ambiguity and perform optimization.
A previous work \cite{zhuKernelMeanEmbedding2020a} already contains the idea of using embedding as optimization variables. However, that work used regularized formulations without rigorously characterizing the distributions corresponding to the solutions.
In kernel methods literature, the work of \cite{kanagawaRecoveringDistributionsGaussiana} used embedding to estimate the probability, but did not solve optimization problems or anticipating the worst case. Their convergence result also requires the function to be in the Besov space, as opposed to Proposition.~\ref{thm:cvgs}.

\section{Kernel mean embedding moment problem}\label{sec:method}
Suppose
$\mu_{\hat P}$ is given by the empirical embedding
$\mu_{\hat P}=\frac1M\sum\phi ( x_i )$
given data samples $\{x_i\}_{i=1}^M$. Alternatively, it can also be the embedding of a continuous distribution such as 
$\mu_{\hat P} = \int  \phi(x)  \ d\hat{P}(x),\ \hat{P}\sim  N(\mu, \sigma)$. 
We now propose
the kernel mean embedding moment problem (KME-MP)
\begin{equation}
        \begin{aligned}&\underset{P,\mu}{\mxz} &  & \int l(x) \ dP(x) \\
                &\sjt& & \|\mu - \mu_{\hat P}\|_\rkhs\leq \epsilon\\
                &&&     \int \phi(x) \ dP(x) = \mu\\
                &&& P\in \mathcal P, \mu\in \rkhs,
        \end{aligned}
\label{eq:kme_mp}
\end{equation}
where 
$l(x)$ is a bounded cost function of interest.
The optimal value of \eqref{eq:kme_mp} is referred to as the worst-case risk.
$\|\cdot\|_\rkhs$ is the Hilbert space norm.
The second constraint is the embedding mapping.
$\mathcal P$ is the set of all probability measures.
Our formulation above makes obvious that, intuitively, the embedding \(\mu\) may be viewed as a \emph{generalized moment vector}. 
Using the \emph{kernel trick}, the RKHS structure allows us to work with infinite-dimensional moment vector \(\mu\),
therefore generalizing the moment problem~\eqref{eq:moment}.


In the optimization problems of this paper, we focus on empirical embeddings of the form $\mu = \kme{\alpha_i}{ z_i}$ with the coefficients $\alpha_i\in\mathbb R$ and expansion points $z_i\in\mathcal X$ (also called locations), for $i=1\dots N$. This can also be seen as a finite-sample estimator of \eqref{eq:defkme}. 
Using this representation, a concrete version of problem~\ref{eq:kme_mp} is given by
\[
\begin{aligned}
&\mxz_{\alpha} & &  \sum_{i=1}^N\alpha_i l(z_i) \\
&\sjt& & \| \kme{\alpha_i}{ z_i} -\sum_{i=1}^M{\frac1M}\phi({x_i})\|_\rkhs\leq\epsilon,
\end{aligned}
\]
where $\{x_i\}_{i=1}^M$ is the data set we have in hand.
Using the kernel trick, the constraint can be written as the convex quadratic constraint
\[
\alpha^\top K_z\alpha - 2\frac1M \alpha^\top K_{zx} \mathbf{1} + \frac{1}{M^2}\mathbf{1}^\top K_x\mathbf{1} \leq\epsilon^2,
\]
where $K_z, K_{zx}, K_x$ are the Gram matrices associated with the kernel $K_z:= [k(z_i, z_j)]_{ij}, K_{zx}:= [k(z_i, x_j)]_{ij}, K_x:= [k(x_i, x_j)]_{ij}$, $\mathbf{1}:=[1,\dots,1]^\top$, and $\alpha := (\alpha_1, \ldots, \alpha_N)^\top$.

\begin{remark}
One may alternatively choose to make $\seqin{z_i}$ decision variables of the optimization formulation, e.g., in machine learning applications \cite{bachEquivalenceHerdingConditional2012}. However, doing so will make the optimization intractable.
\end{remark}


We must first establish a few theoretical results whose detailed proofs are 
presented in Sec.~\ref{sec:theory}.
They lead to the main optimization formulation of this paper. 
Following that, we apply our method to characterize the worst-case constraint violation probability.
\begin{lemma}
\label{thm:kmeisdirac}
Suppose we have an embedding of the form $\mu_P =
\sum_{i=1}^{N}{\alpha_i}\phi({x_i})$ with $N<\infty$, where $\phi$ is a feature map associated with a characteristic kernel.
Then the corresponding signed measure $P$ is a
probability measure if and only if $\sum_{i=1}^N\alpha_i =1,
\alpha_i\geq0$. 
\end{lemma}
Thus, we obtain the practical form of the KME-MP problem
\begin{equation}        
\begin{aligned}
&\mxz_{\alpha} & &  \sum_{i=1}^N\alpha_i l(z_i) \\
&\sjt& &\alpha^\top K_z\alpha - 2\frac1M \alpha^\top K_{zx} \mathbf{1} + \frac{1}{M^2}\mathbf{1}^\top K_x\mathbf{1} \leq\epsilon^2  \\
&&& \sum_{i=1}^N\alpha_i =1, \alpha_i\geq0, i=1\dots N.
\end{aligned}
\label{eq:kmemp_prac}
\end{equation}
This optimization problem is by itself convex and tractable. It does not need a dual reformulation for the purpose of this paper. 
The feasibility of this problem is trivially assured.
\begin{lemma}
        \label{thm:feasible}
        If the new expansion points contain the original data points $\seqin{z_i}\supseteq\{x_i\}_{i=1}^M$ and
        $\epsilon> 0$, the optimization problem \eqref{eq:kmemp_prac} is strictly feasible.
\end{lemma}
Moreover, we have the following relationship between the practical formulation~\eqref{eq:kmemp_prac} and the original formulation~\eqref{eq:kme_mp}.
\begin{prop}
\label{thm:cvgs}
Suppose $l(x)$ is a bounded, lower semicontinuous function.
Let $P^*$ be the solution of the original formulation~\eqref{eq:kme_mp}.
$\{\alpha_i^*\}_{i=1}^N$ is the solution of \eqref{eq:kmemp_prac}
with $\seqin{z_i}$ sampled from a distribution that has common support with $P^*$,
i.e., 
$\seqin{z_i}\sim Q_z,  \mathrm{supp}(Q_z) \supseteq \mathrm{supp}(P^*) $.
Then
$$
\sum_{i=1}^{N}{\alpha_i^*}{l(z_i)} \xrightarrow{N\to\infty} \int l(x) \ dP^*(x).
$$
Furthermore, 
$
\sum_{i=1}^N\alpha_i^* l(z_i)  \leq \int l(x) dP^*(x)
$ for any $N\in\mathbb{Z}^+$,
i.e., the optimal value of \eqref{eq:kmemp_prac} converges to that of \eqref{eq:kme_mp} from below.
\end{prop}
We give a novel proof in Sec.~\ref{sec:proof_prop}. It is important not to misunderstand that $N$ denotes the number of expansion points, not the number of data samples, which is $M$.
{\remark 
In practice, we may not know what $supp(P^*)$ is in order to sample $\seqin{z_i}\sim Q_z$ with common support. 
However, we typically have domain knowledge of the bounds for such samples. For example, we typically know bounds for states of a dynamical system.
This motivates us, in finite-sample cases, to restrict the sampling on a compact domain.
We may simply use a uniform sampling or gridding on some set that contains the whole feasible set. Importantly, we also want to sample points that may violate the constraints, as will be made clear in the next application.
}

This convergence result is particularly useful, e.g. in stress testing models, as we may make the following intuitive statement.
``The worst-case cost is at least as bad as \(\sum_{i=1}^N\alpha_i^* l(z_i)\).''
This is sometimes referred to as an \emph{optimistic} bound.
Hence, our optimization formulation could be understood as searching for a maximized lower bound to the worst-case risk.

\subsection*{Application: evaluating the worst-case probability}\label{sec:chance}
One application of the proposed method is to evaluate the worst-case probability of a certain event. For example, this arises in stochastic control when dealing with a chance constraint of form \(P(X\in C)\geq1-\delta\). In the context of this paper, we consider the presence of ambiguity in the underlying distribution \(P\), i.e., $P\in\mathcal C$ for some ambiguity set $\mathcal C$.
We are interested in evaluating
\[
\begin{aligned}
\mxz_P P(X\notin C)\ \ \sjt\ P\in\mathcal C\subseteq\mathcal P. 
\end{aligned}
\]
Chance constraints with ambiguity is also referred to as distributionally robust chance constraints. Intuitively, by solving this optimization, we are evaluating the worst-case tail probability of the risk.
Obviously, optimization with respect to probability distribution is a special case of the moment problem \eqref{eq:kme_mp} with the cost function \(l(x)\) given by the indicator function \(\indfn\{x\notin C\}\).
Hence, we have the following concrete formulation adapted from \eqref{eq:kmemp_prac}.
\begin{equation}        
\begin{aligned}
&\mxz _ {\alpha} & &  \sum_{i=1}^N \alpha_i \indfn\{z_i\notin C\}\\
&\sjt& &\alpha^\top K_z\alpha - 2\frac1M \alpha^\top K_{zx} \mathbf{1} + \frac{1}{M^2}\mathbf{1}^\top K_x\mathbf{1} \leq\epsilon^2 \\
&&& \sum_{i=1}^N\alpha_i =1, \alpha_i\geq0, i=1\dots N .
\end{aligned}
\label{eq:kme_indicator}
\end{equation}
Our convergence result in Proposition~\ref{thm:cvgs} applies since \(\indfn\{x\notin C\}\) is bounded and lower semicontinuous if $C$ is closed. 
\begin{remark}
(Geometric interpretation)
We may also view the KME expansion $\kme{\alpha_i}{z_i}$ as putting probability mass $\alpha_i$ on location $z_i$. 
Then the worst-case risk can be interpreted as an adversarial, such as nature, trying to move the mass to incur the worst cost. This draws a connection between our approach and optimal transportation metric-based approach such as \cite{mohajerinesfahaniDatadrivenDistributionallyRobust2018b}. We later demonstrate this interpretation in Sec.~\ref{sec:experiment}, Fig.~\ref{fig:transport}.
\end{remark}
\subsection*{Alternative ambiguity set with known first $p$ moments}\label{sec:relation_kmoments}
Ambiguity sets with known first $p$ moments are typically studied separately from metric-ball constraints. 
We now propose the following unified view using KME. 
Let us consider the polynomial kernel of order $p$: \(k(x, x') = (x^\top x' +1) ^p\). For simplicity, we set \(p=2\).
We consider \eqref{eq:kmemp_prac} associated with this  kernel, obtaining the following moment problem under known first $p$ moments. 
\begin{equation}
\begin{aligned}
&\mxz _ {P\in {\mathcal P},\mu_P} &  &\sum_{i=1}^N\alpha_i l(z_i) 
\\&\sjt& &\alpha^\top K_z\alpha \\
&&&- 2 \sum_{i=1}^N\alpha_i({z_i}^\top \mathbb{E}{x x^\top} z_i + 2 \mathbb{E} {x}^\top z_i +1) \\
&&& + \mathrm{Tr}(\mathbb{E}{x x^\top}\mathbb{E}{x x^\top}  ) + 2 \mathbb{E} {x}^\top\mathbb{E} {x} +1
=0
\\
&&& \sum_{i=1}^N\alpha_i =1, \alpha_i\geq0, i=1\dots N,
\end{aligned}
\label{eq:mp_pol}
\end{equation}
where $K_z:= [k(z_i, z_j)]_{ij}$ is the Gram matrix.
However, the polynomial kernel is not characteristic: the mapping $\int \phi(x) \ dP(x) = \mu_P$ is finite-dimensional. Intuitively speaking, it can not distinguish all distributions that share the (finitely many) known moments, hence incurring ambiguity. Problem~\eqref{eq:mp_pol} is convex, nonetheless.
\begin{remark}
The first $p$ moments are typically estimated from data.
Therefore, it is sometimes argued that one may skip the estimation step and favor directly using nonparametric formulation such as KME-MP~\eqref{eq:kmemp_prac}. 
For this reason as well as clarity, this paper focuses on \eqref{eq:kmemp_prac} instead of \eqref{eq:mp_pol}.
\end{remark}

\section{Numerical experiment}\label{sec:experiment}
In this section,
we solve the proposed optimization problem \eqref{eq:kmemp_prac}, \eqref{eq:kme_indicator} to quantify the worst-case risk in synthetic example and a stochastic control problem.
We focus on describing experimental results and their implications. Additional detailed descriptions of the experimental setup are given in the
appendix.
\subsection{Synthetic worst-case probability problem}
\label{sec:toy}
Suppose $X$ is a real-valued random variable whose distribution is \emph{unknown} to us.
We only have access to the empirical data samples $\{x_i\}_{i=1}^M$.
We are interested in quantifying the worst-case violation probability of the event $\{X\leq c\}$.
To this end, we cast the problem in the form of \eqref{eq:kme_mp}.
$$
\begin{aligned}
&\mxz _ {P\in\mathcal P, \mu\in\rkhs} & &  P(X > c)\\
&\sjt& &\|\mu-\mu_{\hat P}\|_\rkhs\leq\epsilon, 
\int \phi(x) \ dP(x) = \mu,
\end{aligned}
$$
where $\mu_{\hat P}=\sum_{i=1}^M\frac1M\phi(x_i)$ is the mean embedding of the empirical data distribution.

\begin{figure}[tb!]
        \centering
                \includegraphics[width=0.75\columnwidth]{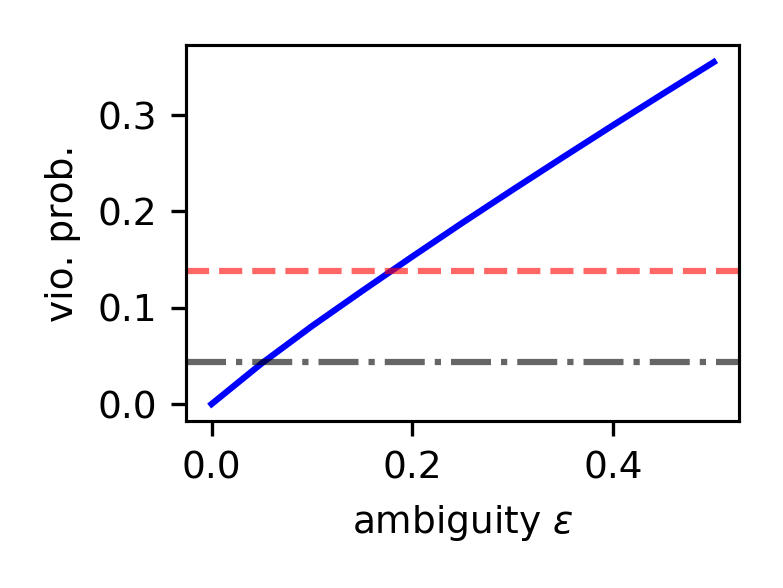}
                \caption{(blue line) This figure depicts the value of worst-case violation probability ($y$-axis), i.e., \(\sup_P P(x> c)\) changes with respect to the various allowed ambiguity levels \(\epsilon\) ($x$-axis) in the worst-case risk quantification problem. As we increase the ambiguity level $\epsilon$, the worst-case violation probability exceeds both Chernoff's (gray dashed) and Cantelli's (red dashed) bound.
                }
                \label{fig:cantelli_n100}
\end{figure}

To estimate the worst-case violation probability, we solve the proposed optimization formulation~\eqref{eq:kme_indicator} with various levels of ambiguity $\epsilon$. 
As a result, we obtain the solution to \eqref{eq:kme_indicator}, $\seqin{\alpha_i^*}$, as well as the corresponding worst-case violation probability.
We plot this quantity in Fig.~\ref{fig:cantelli_n100}. 
As we can observe from the figure, the proposed method is able to estimate the worst-case value as the ambiguity level $\epsilon$ increases.
See the figure caption for more details.

Let us first recall classical Cantelli's inequality mentioned in the introduction.
It bounds the tail probability of any distribution sharing the first two moments. 
In our dataset, Cantelli's inequality yields a tail bound
$
P(x\geq 2.5)\leq\frac{1}{1+2.5^2} \approx 0.138,
$
which is of course not sharp.
For example, a sharper bound in this case is given by the Chernoff's bound.
$
P(x\geq 2.5)\leq \exp(-2.5^2/2)\approx 0.044.
$
Even pessimistic bounds are subject to violation under distributional ambiguity. As we can observe in Fig.~\ref{fig:cantelli_n100}, with a large enough ambiguity level \(\epsilon\), the worst-case risk exceeds that of the bounds given by the Chernoff's and Cantelli's inequality. This implies the importance of anticipating distributional ambiguity that may defeat even conservative safety bounds.
\wjsay{I am actually confused about the main point of comparing the result to
these bounds. Is it correct that for these bounds, you use
$\mathcal{N}(0,1)$? Exceeding these bounds is not surprising to me at all
since you consider other distributions around $\mathcal{N}(0,1)$ as well.}
\jz{True. Need to present this better. My main point is that, in reality, we do not know the distribution $N(0,1)$, but only the data. Then if there's distribution shift, we can be quite far from any safety bounds.}

\begin{figure}[tb!]
        \centering
        \includegraphics[width=0.75\columnwidth]{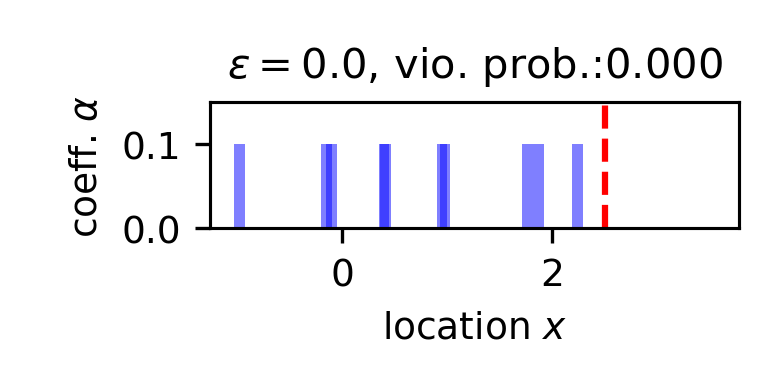}
        \includegraphics[width=0.75\columnwidth]{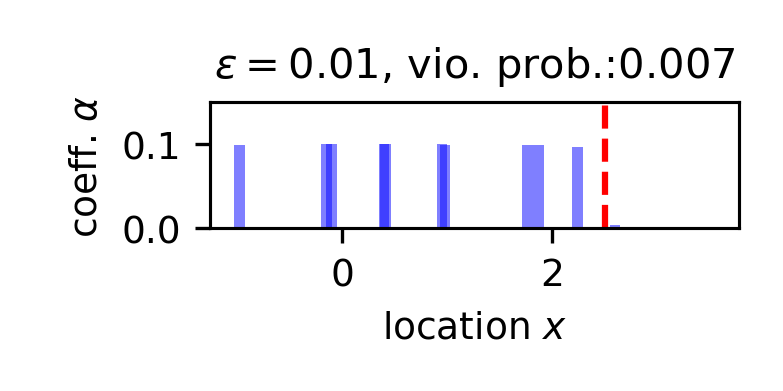}
        \includegraphics[width=0.75\columnwidth]{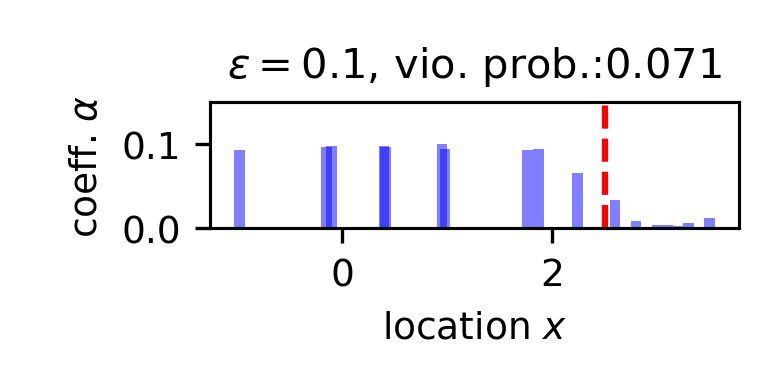}
        \caption{Examining the worst-case constraint violation behavior under different ambiguity levels. The (blue) vertical bars indicate the mass (height) allocated at a certain location $x$ (horizontal axis). The (red) bar indicates the constraint level $c=2.5$. (top) $\epsilon=0$, no distributional ambiguity is allowed. All the mass stayed within the constraint. (middle) $\epsilon=0.01$, a small amount of mass is shifted to the right side. The violation probability is $0.7\%$. (bottom) $\epsilon=0.1$, violation probability $7.1\%$. More mass is allowed to be moved to the violation side due to a larger tolerance of ambiguity.
        }
        \label{fig:transport}
\end{figure}
We now analyze the behavior of our approach more carefully. 
As we have discussed in Sec.~\ref{sec:method}, we may interpret the coefficient $\alpha_i$ of the embedding $\kme{\alpha_i}{z_i}$ as the mass of this discrete distribution at location $x_i$.
For simplicity of viewing, we first sample a small set of \(M=10\) points \(\{x_i\}_{i=1}^M\) to be our empirical data samples. In this case, all sampled points are on the left side of \(c=2.5\). The naive Monte Carlo estimate of the violation probability yields
$ P(X > c)\approx\sum_{i=1}^N\frac1N\indfn\{x_i> c\} = 0$.
Recall this is the phenomenon of silent evidence discussed in the introduction.
We now demonstrate that, by applying the proposed method, we are able to anticipate the potential violation with limited data.

We then solve the optimization problem~\eqref{eq:kme_indicator} 
and plot the expansion weights \(\seqin{\alpha_i}\) in Fig.~\ref{fig:transport} for various ambiguity levels $\epsilon$.
Similar to our geometric interpretation of \eqref{eq:kme_indicator}, the mass is ``\emph{transported}'' from left-hand-side of \(c\) (constraint satisfaction), to the right-hand-side (constraint violation). This is reminiscent of the worst-case risk characterized by optimal transport metric, e.g., in \cite{mohajerinesfahaniDatadrivenDistributionallyRobust2018b}.
\begin{figure}[tb!]
        \centering
        \includegraphics[width=0.75\columnwidth]{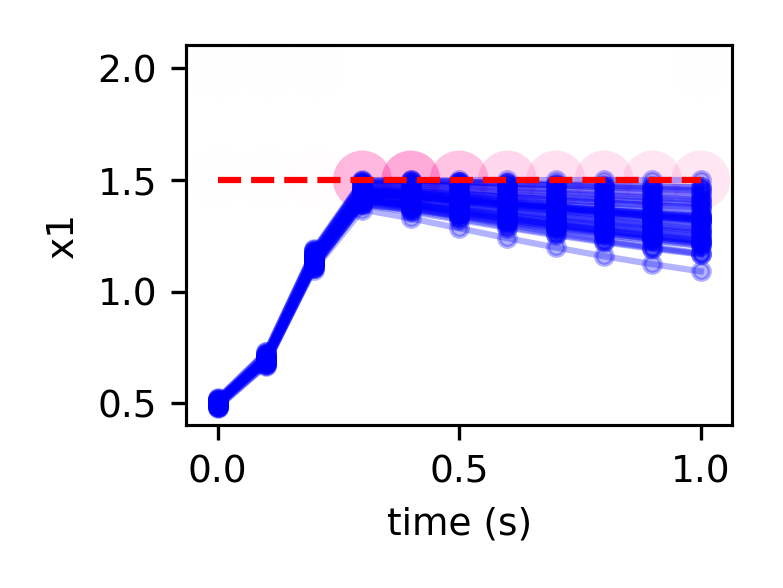}
        \includegraphics[width=0.75\columnwidth]{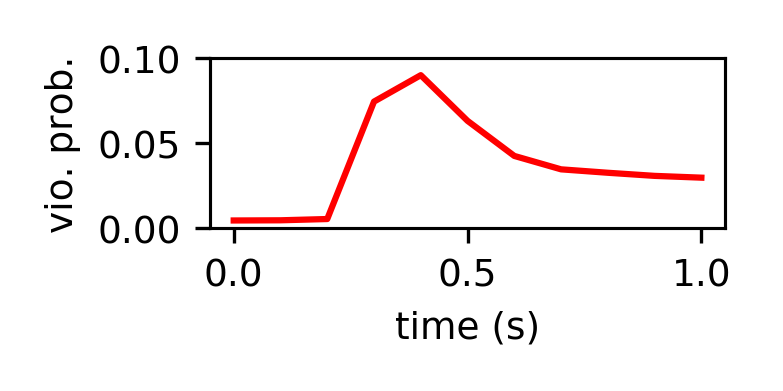}
        \caption{This figure is best viewed in color. In this simulation, the ambiguity level is set to $\epsilon=0.01$. (top) State trajectory (blue) of the solution to the stochastic OCP. The constraint is to stay below the (red dashed) line. The (pink) shaded area denotes violation mass is transported to this location. The higher the mass, the darker the shade.
        (bottom) is the plot of the worst-case violation probability $\sup_P P(x>1.5)$ under ambiguity by solving the optimization problem~\eqref{eq:kme_indicator}. The horizontal axis is time. 
        }
        \label{fig:smpc}
\end{figure}
\subsection{Worst-case constraint violation in stochastic control}
\label{sec:smpc}
In this experiment, we are interested in the constrained stochastic optimal control problem (OCP). We consider the system model of the Van der Pol oscillator
\[
\frac{d}{dt}{\left[\begin{array}{l}   {x_{1}} \\  {x_{2}} \end{array}\right]}=\left[\begin{array}{c} {x_{2}} \\  {-0.1\left(1-x_{1}^{2}\right) x_{2}-x_{1}+u}    \end{array}\right]. \label{eq:vanderpol}
\]
This example was adapted from \cite{anderssonCasADiSoftwareFramework2019}. The goal of the control design is to steer the system state \(x_1\) as high as possible while staying below the bound (for simplicity, we use the constraint \(x_{1} \leq 1.5\)). This is formulated as the following uncertain OCP.
\begin{equation}
\begin{array}{ll}   {\ \underset{x(\cdot), u(\cdot)}{\operatorname{minimize}}} & {\int_{0}^{T}\left\|x_{1}(t)-3.0\right\|_{2}^{2} d t} \\
         {\text { subject to }} & {\dot{x}(t)=f(x(t), u(t)), \quad \forall t \in[0, T]} \\    {} & {-40 \leq u(t) \leq 40, \quad \forall t \in[0, T]} \\   {} & -0.25 \leq x_{1}(t)\leq 1.5,\quad  \forall t \in[0, T]\\ {} & {x(0)=s},   \end{array} 
\label{eq:ocp}
\end{equation}
where the initial state $s$ is uncertain.
To work with the nonlinear dynamics, we adopt the scenario MPC~\cite{calafioreScenarioApproachRobust2006} to solve for the optimal control input.
The state trajectory associated with the OCP solution is plotted in Fig.~\ref{fig:smpc} (top).

Furthermore, we assume there is ambiguity in the state distribution \(\mathbb P_{x_t}\). Given the empirical data distribution \(\hat{\mathbb P}_{x_t}=\frac1N\sum\delta_{{x_t}^{(i)}}\) at time \(t\), we are interested in quantifying the worst-case violation probability
\[
\sup_P P(X_t > 1.5)\ \sjt \|\mu_P - \mu_{\hat{\mathbb P}_{x_t}}\|_{\mathcal H}\leq\epsilon.
\]
The recent work of \cite{zhuNewDistributionFreeConcept2019a} has advocated using embedding for representing and propagating uncertainty. We now demonstrate our framework can use this representation to stress test the control design.
To this end, we solve \eqref{eq:kme_indicator} and plot the worst-case violation probability at each time-step in Fig.~\ref{fig:smpc} (bottom). 
There is no violation in the beginning as the initial states are far from the constraint boundary. However, as the system is steered closer to the boundary, the worst-case violation probability rises.

Similar to the previous example, we visualize the geometric interpretation in Fig.~\ref{fig:smpc} (top), illustrated as the shaded area (pink).
The darkness of the shade is proportional to the amount of the violation probability mass.
We may interpret our method as transporting mass to the locations of violation states to create a worst-case distribution. 
\begin{figure}[tb!]
        \centering
        \includegraphics[width=0.75\columnwidth]{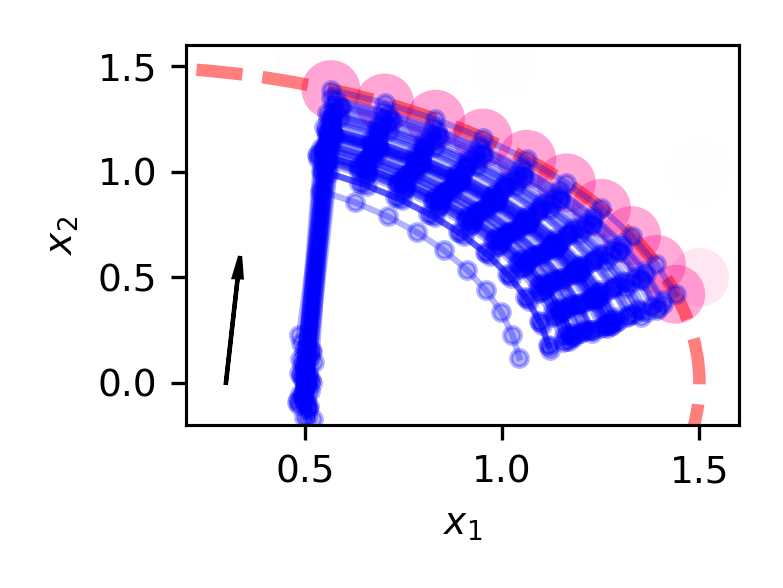}
        \includegraphics[width=0.75\columnwidth]{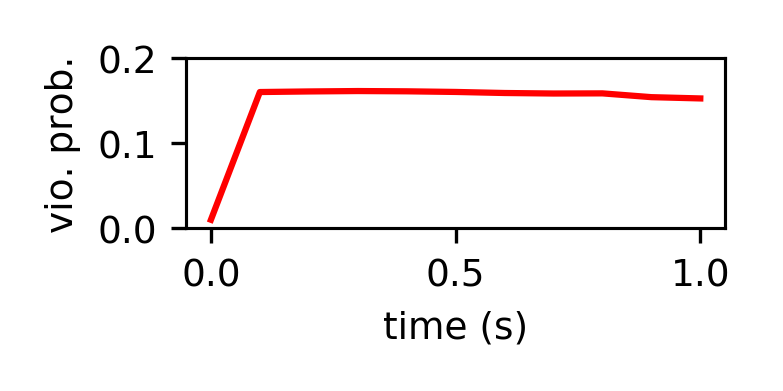}
        \caption{This figure is best viewed in color. 
        The dynamics is the same as that in Fig.~\ref{fig:smpc} except for we try to quantify the worst-case violation probability for the non-polyhedral constraint regions, $\sup_P P\{\sqrt{(x_1^2+x_2^2)}> 1.5\}$ under ambiguity. (top) is the phase plot of the system. (arrow) indicates the direction of the state trajectory over time. The (red) boundary is a circle. 
        The (black) arrow indicates the direction of the state over time.
        The worst-case violation is shaded in (pink) color. (bottom) plots the violation probability in time.
        }
        \label{fig:phase}
\end{figure}
Unlike previous approaches such as \cite{vanparysGeneralizedGaussInequalities2016}, our method does not require the constraint region to be polyhedral. We now test the case where the event has a simple non-nonlinear multivariate constraint $\{\sqrt{(x_1^2+x_2^2)}\leq 1.5\}$, where $x_1, x_2$ are the two coordinates of the state. (In fact, it does not need to be quadratic either.) Fig.~\ref{fig:phase} illustrates the resulting phase plot of the system and the worst-case violation probability. The interpretation is similar to Fig.~\ref{fig:smpc}.

The implication of this experiment is interesting from the perspective of optimal control design. 
Traditional stochastic control aims to bound most of the (known) probability mass of the uncertainty to stay within bound --- such as constraint tightening in stochastic MPC (cf. \cite{farinaStochasticLinearModel2016,mesbahStochasticModelPredictive2016a}). 
We aim to go beyond such constraint tightening methods based on known probability measures.
Although all the stochastic control trajectories in our experiment satisfy the constraints, we were able to quantify non-zero worst-case violation probability.
To robustify against the worst-case risk under the distributional ambiguity quantified by this paper, the control design must be \emph{distributionally robust}.

\section{Conclusion}\label{sec:conclusion}
In this paper, we introduced the kernel mean embedding framework as a tool for solving generalized moment problems. 
We presented constrained optimization formulations with the ambiguity sets described by RKHS distance as well as practical solution methods.
Through theoretical analysis and numerical experiments,
we demonstrate how the proposed method can characterize the worst-case risk under distributional ambiguity.
Therefore, we view the kernel mean embedding as a promising tool for robust optimization and control under distributional ambiguity.
One future direction of work is
choosing the size $\epsilon$ of the ambiguity set, which can be motivated by
nonparametric statistical tests~\cite{grettonKernelTwosampleTest2012}.
Another direction is exploring the connection of our approach with the popular Wasserstein metric, possibly by using the \emph{integral probability metric} framework.
\section{Acknowledgment}\label{acknowledgment}

We thank the anonymous reviewers for their constructive suggestions during the review process. 

\bibliographystyle{IEEEtran}
\bibliography{cdc_uq}

\begin{thebibliography}{10}
\providecommand{\url}[1]{#1}
\csname url@rmstyle\endcsname
\providecommand{\newblock}{\relax}
\providecommand{\bibinfo}[2]{#2}
\providecommand\BIBentrySTDinterwordspacing{\spaceskip=0pt\relax}
\providecommand\BIBentryALTinterwordstretchfactor{4}
\providecommand\BIBentryALTinterwordspacing{\spaceskip=\fontdimen2\font plus
\BIBentryALTinterwordstretchfactor\fontdimen3\font minus
  \fontdimen4\font\relax}
\providecommand\BIBforeignlanguage[2]{{%
\expandafter\ifx\csname l@#1\endcsname\relax
\typeout{** WARNING: IEEEtran.bst: No hyphenation pattern has been}%
\typeout{** loaded for the language `#1'. Using the pattern for}%
\typeout{** the default language instead.}%
\else
\language=\csname l@#1\endcsname
\fi
#2}}

\bibitem{talebBlackSwanImpact2007}
N.~N. Taleb, \emph{The Black Swan: {{The}} Impact of the Highly
  Improbable}.\hskip 1em plus 0.5em minus 0.4em\relax {Random house}, 2007,
  vol.~2.

\bibitem{popescuSemidefiniteProgrammingApproach2005}
I.~Popescu, ``\BIBforeignlanguage{en}{A {{Semidefinite Programming Approach}}
  to {{Optimal}}-{{Moment Bounds}} for {{Convex Classes}} of
  {{Distributions}}},'' \emph{\BIBforeignlanguage{en}{Mathematics of Operations
  Research}}, vol.~30, no.~3, pp. 632--657, Aug. 2005.

\bibitem{farinaStochasticLinearModel2016}
M.~Farina, L.~Giulioni, and R.~Scattolini, ``\BIBforeignlanguage{en}{Stochastic
  linear {{Model Predictive Control}} with chance constraints \textendash{}
  {{A}} review},'' \emph{\BIBforeignlanguage{en}{Journal of Process Control}},
  vol.~44, pp. 53--67, Aug. 2016.

\bibitem{sinhaCertifyingDistributionalRobustness2017a}
A.~Sinha, H.~Namkoong, and J.~Duchi, ``Certifying {{Some Distributional
  Robustness}} with {{Principled Adversarial Training}},'' pp. 1--41, 2017.

\bibitem{isiiSharpnessTchebychefftypeInequalities1962}
K.~Isii, ``\BIBforeignlanguage{en}{On sharpness of tchebycheff-type
  inequalities},'' \emph{\BIBforeignlanguage{en}{Annals of the Institute of
  Statistical Mathematics}}, vol.~14, no.~1, pp. 185--197, Dec. 1962.

\bibitem{shapiroDualityTheoryConic2001}
A.~Shapiro, ``\BIBforeignlanguage{en}{On {{Duality Theory}} of {{Conic Linear
  Problems}}},'' in \emph{\BIBforeignlanguage{en}{Semi-{{Infinite
  Programming}}}}, P.~Pardalos, M.~{\'A}. Goberna, and M.~A. L{\'o}pez,
  Eds.\hskip 1em plus 0.5em minus 0.4em\relax {Boston, MA}: {Springer US},
  2001, vol.~57, pp. 135--165.

\bibitem{rogosinskiMomentsNonNegativeMass}
W.~W. Rogosinski, ``\BIBforeignlanguage{en}{Moments of {{Non}}-{{Negative
  Mass}}},'' p.~28.

\bibitem{lasserreBoundsMeasuresSatisfying2002}
J.~B. Lasserre, ``Bounds on measures satisfying moment conditions,'' \emph{The
  Annals of Applied Probability}, vol.~12, no.~3, pp. 1114--1137, 2002.

\bibitem{bertsimasMomentProblemsSemidefinite2000}
D.~Bertsimas and J.~Sethuraman, ``Moment problems and semidefinite
  optimization,'' in \emph{Handbook of Semidefinite Programming}.\hskip 1em
  plus 0.5em minus 0.4em\relax {Springer}, 2000, pp. 469--509.

\bibitem{bertsimasOptimalInequalitiesProbability2005}
D.~Bertsimas and I.~Popescu, ``\BIBforeignlanguage{en}{Optimal {{Inequalities}}
  in {{Probability Theory}}: {{A Convex Optimization Approach}}},''
  \emph{\BIBforeignlanguage{en}{SIAM Journal on Optimization}}, vol.~15, no.~3,
  pp. 780--804, Jan. 2005.

\bibitem{vandenbergheGeneralizedChebyshevBounds2007}
L.~Vandenberghe, S.~Boyd, and K.~Comanor, ``Generalized {{Chebyshev Bounds}}
  via {{Semidefinite Programming}},'' \emph{SIAM Review}, vol.~49, no.~1, pp.
  52--64, Jan. 2007.

\bibitem{vanparysGeneralizedGaussInequalities2016}
B.~P.~G. Van~Parys, P.~J. Goulart, and D.~Kuhn,
  ``\BIBforeignlanguage{en}{Generalized {{Gauss}} inequalities via semidefinite
  programming},'' \emph{\BIBforeignlanguage{en}{Mathematical Programming}},
  vol. 156, no. 1-2, pp. 271--302, Mar. 2016.

\bibitem{delageDistributionallyRobustOptimization2010}
E.~Delage and Y.~Ye, ``\BIBforeignlanguage{en}{Distributionally {{Robust
  Optimization Under Moment Uncertainty}} with {{Application}} to
  {{Data}}-{{Driven Problems}}},'' \emph{\BIBforeignlanguage{en}{Operations
  Research}}, vol.~58, no.~3, pp. 595--612, June 2010.

\bibitem{erdoganAmbiguousChanceConstrained2006}
E.~Erdo{\u g}an and G.~Iyengar, ``\BIBforeignlanguage{en}{Ambiguous chance
  constrained problems and robust optimization},''
  \emph{\BIBforeignlanguage{en}{Mathematical Programming}}, vol. 107, no.~1,
  pp. 37--61, June 2006.

\bibitem{mohajerinesfahaniDatadrivenDistributionallyRobust2018b}
P.~Mohajerin~Esfahani and D.~Kuhn, ``\BIBforeignlanguage{en}{Data-driven
  distributionally robust optimization using the {{Wasserstein}} metric:
  Performance guarantees and tractable reformulations},''
  \emph{\BIBforeignlanguage{en}{Mathematical Programming}}, vol. 171, no.~1,
  pp. 115--166, Sept. 2018.

\bibitem{ben-talRobustSolutionsOptimization2013}
A.~{Ben-Tal}, D.~{den Hertog}, A.~De~Waegenaere, B.~Melenberg, and G.~Rennen,
  ``\BIBforeignlanguage{en}{Robust {{Solutions}} of {{Optimization Problems
  Affected}} by {{Uncertain Probabilities}}},''
  \emph{\BIBforeignlanguage{en}{Management Science}}, vol.~59, no.~2, pp.
  341--357, Feb. 2013.

\bibitem{bertsimasDatadrivenRobustOptimization2017}
D.~Bertsimas, N.~Kallus, and V.~Gupta, \emph{Data-Driven Robust
  Optimization}.\hskip 1em plus 0.5em minus 0.4em\relax {Springer Berlin
  Heidelberg}, 2017, publication Title: Mathematical Programming.

\bibitem{owhadi2013optimal}
H.~Owhadi, C.~Scovel, T.~J. Sullivan, M.~McKerns, and M.~Ortiz, ``Optimal
  uncertainty quantification,'' \emph{Siam Review}, vol.~55, no.~2, pp.
  271--345, 2013.

\bibitem{berlinetReproducingKernelHilbert2011}
A.~Berlinet and C.~{Thomas-Agnan}, \emph{Reproducing Kernel {{Hilbert}} Spaces
  in Probability and Statistics}.\hskip 1em plus 0.5em minus 0.4em\relax
  {Springer Science \& Business Media}, 2011.

\bibitem{smolaHilbertSpaceEmbedding2007}
A.~Smola, A.~Gretton, L.~Song, and B.~Sch{\"o}lkopf, ``A hilbert space
  embedding for distributions,'' \emph{Lecture Notes in Computer Science
  (including subseries Lecture Notes in Artificial Intelligence and Lecture
  Notes in Bioinformatics)}, vol. 4754 LNAI, pp. 13--31, 2007.

\bibitem{sriperumbudurUniversalityCharacteristicKernels2011}
B.~K. Sriperumbudur, K.~Fukumizu, and G.~R.~G. Lanckriet, ``Universality,
  {{Characteristic Kernels}} and {{RKHS Embedding}} of {{Measures}},''
  \emph{Journal of Machine Learning Research}, vol.~12, no. Jul, pp.
  2389--2410, 2011.

\bibitem{grettonKernelTwosampleTest2012}
A.~Gretton, K.~M. Borgwardt, M.~J. Rasch, B.~Sch{\"o}lkopf, and A.~Smola, ``A
  kernel two-sample test,'' \emph{Journal of Machine Learning Research},
  vol.~13, no. Mar, pp. 723--773, 2012.

\bibitem{zhuKernelMeanEmbedding2020a}
J.-J. Zhu, B.~Sch{\"o}lkopf, and M.~Diehl, ``A {{Kernel Mean Embedding
  Approach}} to {{Reducing Conservativeness}} in {{Stochastic Programming}} and
  {{Control}},'' \emph{arXiv preprint arXiv:2001.10398}, 2020.

\bibitem{kanagawaRecoveringDistributionsGaussiana}
M.~Kanagawa and K.~Fukumizu, ``\BIBforeignlanguage{en}{Recovering
  {{Distributions}} from {{Gaussian RKHS Embeddings}}},'' p.~9.

\bibitem{bachEquivalenceHerdingConditional2012}
F.~Bach, S.~{Lacoste-Julien}, and G.~Obozinski, ``On the {{Equivalence}}
  between {{Herding}} and {{Conditional Gradient Algorithms}},''
  \emph{arXiv:1203.4523 [cs, math, stat]}, Sept. 2012.

\bibitem{anderssonCasADiSoftwareFramework2019}
J.~A.~E. Andersson, J.~Gillis, G.~Horn, J.~B. Rawlings, and M.~Diehl,
  ``\BIBforeignlanguage{en}{{{CasADi}}: A software framework for nonlinear
  optimization and optimal control},''
  \emph{\BIBforeignlanguage{en}{Mathematical Programming Computation}},
  vol.~11, no.~1, pp. 1--36, Mar. 2019.

\bibitem{calafioreScenarioApproachRobust2006}
G.~Calafiore and M.~Campi, ``The scenario approach to robust control design,''
  \emph{IEEE Transactions on Automatic Control}, vol.~51, no.~5, pp. 742--753,
  May 2006, conference Name: IEEE Transactions on Automatic Control.

\bibitem{zhuNewDistributionFreeConcept2019a}
J.-J. Zhu, K.~Muandet, M.~Diehl, and B.~Sch{\"o}lkopf, ``A {{New
  Distribution}}-{{Free Concept}} for {{Representing}}, {{Comparing}}, and
  {{Propagating Uncertainty}} in {{Dynamical Systems}} with {{Kernel
  Probabilistic Programming}},'' \emph{arXiv:1911.11082 [cs, eess, math,
  stat]}, Nov. 2019.

\bibitem{mesbahStochasticModelPredictive2016a}
A.~Mesbah, ``Stochastic {{Model Predictive Control}}: {{An Overview}} and
  {{Perspectives}} for {{Future Research}},'' \emph{IEEE Control Systems
  Magazine}, vol.~36, no.~6, pp. 30--44, Dec. 2016, conference Name: IEEE
  Control Systems Magazine.

\bibitem{wachterImplementationPrimalDual}
A.~Wachter and L.~T. Biegler, ``On the implementation of a
  primal\textemdash{}dual interior point filter line search algorithm for
  large-scale nonlinear programming, mathematical programming,'' \emph{Math.
  Program}, vol. 106, no.~1.

\bibitem{pedregosaScikitlearnMachineLearning2011}
F.~Pedregosa, G.~Varoquaux, A.~Gramfort, V.~Michel, B.~Thirion, O.~Grisel,
  M.~Blondel, P.~Prettenhofer, R.~Weiss, and V.~Dubourg, ``Scikit-learn:
  {{Machine}} learning in {{Python}},'' \emph{Journal of machine learning
  research}, vol.~12, no. Oct, pp. 2825--2830, 2011.

\bibitem{parthasarathyProbabilityMeasuresMetric2014}
K.~R. Parthasarathy, \emph{\BIBforeignlanguage{en}{Probability {{Measures}} on
  {{Metric Spaces}}}}.\hskip 1em plus 0.5em minus 0.4em\relax {Academic Press},
  July 2014.

\end{thebibliography}

\section{Appendix: experimental setup}\label{sec:appendix}
All the experiments are implemented using Python with optimization solver
IPOPT~\cite{wachterImplementationPrimalDual}. The stochastic control
experiment used the CasADi
library~\cite{anderssonCasADiSoftwareFramework2019}. 
We used the Gaussian kernel $k(x,y) = \exp\left( -\frac{ \| x-y\|^2_2}{2 \sigma^2}
\right)$  with bandwidth $\sigma$ set to be $\mathrm{median}(\{ \|x_i - x_j \|_2/\sqrt{2} \mid i,j=1,\ldots, N \})$.
This is a common practice known as the median heuristic \cite{grettonKernelTwosampleTest2012}.
We also used sum-of-kernel practice to combine different widths to form a new kernel.
Our specific implementation of kernel computation is adapted from \cite{pedregosaScikitlearnMachineLearning2011}.
\jz{@WJ: could you check this part for kernel width? basically I used median heuristic, and sum of kernels used scaling of the widths, $\sigma\times [0.01, 0.1, 1.0, 10, 100]$} 
\wjsay{The scaling factors appear to be very broad. You might want to make a
bit denser grid. If the data are unimodal, the median heuristic usually works
well.}


\section{Proofs of theoretical results}
\label{sec:theory}
\subsection{Proof of Lemma~\ref{thm:kmeisdirac}}
\begin{proof}
Consider a convex combination of $N$ Dirac measures $P = \sum^{N}_{i=1} \alpha_i \delta_{z_i}$. Its embedding is $\mu_P = \sum_{i=1}^{N}{\alpha_i}\phi({z_i})$. Obviously, $P$ is a probability measure if and only if $\sum_{i=1}^N\alpha_i =1, \alpha_i\geq0$. 
Since the embedding is injective if $\phi$ is associated with a universal kernel,
the conclusion follows.       
\end{proof}

\subsection{Proof of Lemma~\ref{thm:feasible}}
\begin{proof}
        Since $\{x_i\}_{i=1}^M\subseteq\seqin{z_i}$ 
, we can trivially choose $\sum_{i=1}^N\alpha_i\phi(z_i) = \sum_{i=1}^M\frac1M\phi(x_i)$ which is strictly feasible.        
\end{proof}
\subsection{Proof of Proposition~\ref{thm:cvgs}}
\label{sec:proof_prop}
\begin{proof}
\newcommand{\dn}{\mathcal D_N}
\newcommand{\qe}{Q_{n_\epsilon}}
\newcommand{\zseq}{z_i^{(n),N}} 
Let $\mathcal{D}$ be the space of all convex combinations of Dirac measures supported on $\mathcal X$.
It is a fact $\mathcal{D}$ is dense in the space of all probability distributions $\mathcal P$ on $\mathcal X$. Its proof can be found in standard texts, e.g., \cite{parthasarathyProbabilityMeasuresMetric2014} Theorem~{6.3}. Hence, there exists a sequence of probability measures  $Q_n\in \mathcal D$ that converges to $P^*$ weakly, i.e. $Q_n\Rightarrow P^*$. 
Let $Q_n=\sum_{i=1}^{K^{(n)}}{\gamma_i^{(n)}}{\delta_{y_i^{(n)}}}$, where$\sum_{i=1}^{K^{(n)}}{\gamma_i^{(n)}}=1, {\gamma_i^{(n)}}\geq 0$ for $i=1\dots K^{(n)}$. By the Portmanteau theorem, this weak convergence implies
\begin{equation}
    \liminf_{n\to\infty}\int l(x) \ d Q_n(x)
    \geq
    \int l(x) \ dP^*(x),
\label{eq:wk_cvgs}
\end{equation}
for any bounded lower semicontinuous function $l$. 
Note we assume $\mathrm{supp}(Q_n)\subseteq\mathrm{supp}(P^*)$, since otherwise we could replace $Q_n$ by its restriction to $\mathrm{supp}(P^*)$.

Now we prove that $\liminf_{N\to\infty} \sum_{i=1}^N\alpha_i^* l(z_i)\geq\int l(x) \ d Q_n (x)$ for a given $n$.
Due to the assumptions that the $\seqin{z_i}$ are generated from the distribution whose support contains $\mathrm{supp}(P^*)$, 
we have for any support point $y_i^{(n)}$ of $Q_n$,
$$
\mathrm{dist}(\{z_i\}_{i=1}^N, y_i^{(n)}):=
\min_{z\in \{z_i\}_{i=1}^N}d(z,y_i^{(n)}) \xrightarrow{N\to\infty} 0,
$$
where $d$ is the usual distance on $\mathcal X$. Therefore, there exists a sequence $\{\zseq\}_{N=1}^{\infty}$ \wjsay{could you please write which index varies? $i$?}\jz{N varies, i is basically matching that of $y_i$} such that $\zseq\in \seqin{z_i}$ and
\begin{equation}
d(\zseq, y_i^{(n)})\xrightarrow{N\to\infty} 0 .
\label{eq:ztoy}
\end{equation}
For fixed $n$, we use $\{\zseq\}_{N=1}^{\infty}$ to construct a sequence of distributions $\hat
Q_n^N:=\sum_{i=1}^{K^{(n)}}{\gamma_i^{(n)}}\delta_{\zseq}$ that will converge weakly to $Q_n$, where ${K^{(n)}}, {\gamma_i^{(n)}}$ are from the definition of $Q_n$. \wjsay{Sum over what index? From what to what? Also
do the weights sum to one?}
\jz{added to the first time $\gamma$ appears when I introduced $Q_n$}
Due to the lower semicontinuity of $l$ and \eqref{eq:ztoy}, we have
\begin{multline}  
    \liminf_{N\to\infty}\int l(x) \ d\hat Q_n^N(x)
    =
       \sum_{i=1}^{K^{(n)}}{\gamma_i^{(n)}}\left(\liminf_{N\to\infty}l(\zseq)\right) \\
    \geq
        \sum_{i=1}^{K^{(n)}}
        {\gamma_i^{(n)}}
        l (y_i^{(n)})
    =
    \int l(x) \ dQ_n(x).
    \label{eq:qhat}
\end{multline}

As $\seqin{\alpha^*_i}$ attains the maximum of \eqref{eq:kmemp_prac} and $\zseq\in \seqin{z_i}$,
$$
\sum_{i=1}^N\alpha_i^* l(z_i)
\geq\sum_{i=1}^{K^{(n)}}{\gamma_i^{(n)}}{l(\zseq)}.
$$
In the above inequality, we let $N\to\infty$ and use \eqref{eq:qhat},
$$
\liminf_{N\to\infty}
\sum_{i=1}^N\alpha_i^* l(z_i)
\geq\sum_{i=1}^{K^{(n)}}{\gamma_i^{(n)}}{l(y_i^{(n)})}=\int l(x)\ dQ_n(x).
$$

Now we let $n\to\infty$ and use \eqref{eq:wk_cvgs}, we obtain
\begin{multline}
\liminf_{N\to\infty}
\sum_{i=1}^N\alpha_i^* l(z_i)
\geq 
\liminf_{n\to\infty}\int l(x) \ d Q_n(x)\\
\geq
\int l(x)\ dP^*(x).
\label{eq:geqPstar}
\end{multline}
Since $P^*$ is the probability measure that attains the worst-case objective value, we always have 
\begin{equation}
\sum_{i=1}^N\alpha_i^* l(z_i)
\leq 
\int l(x)\ dP^*(x),
\forall N\in\mathbb Z^+.
\label{eq:leqPstar}
\end{equation}
Combining \eqref{eq:geqPstar} and \eqref{eq:leqPstar}, we obtain the convergence.

\begin{todo}
the inequality is in fact equality. Therefore,
$$
\lim_{N\to\infty}\sum_{i=1}^N\alpha_i^* l(z_i) = \int l(x) \ dP^*(x).
$$

Therefore, for any $\epsilon>0$, there exists $\qe=\sum_{i=1}^{K^{(n_\epsilon)}}{\gamma_i^{(n_\epsilon)}}{l({y_i^{(n_\epsilon)})}}\in\mathcal D$ such that 
\begin{equation}
|\int l(x) \ d\qe (x)-\int l(x) \ dP^*(x)|<\epsilon
\label{eq:prop1}
\end{equation}

Now we prove that $\liminf_{N} \sum_{i=1}^N\alpha_i^* l(z_i)\geq\int l(x) \ d\qe (x)$
Let $\dn$ be the set of all distributions that correspond to the feasible set of the optimization problem \eqref{eq:kmemp_prac}. 
It is a discrete distribution supported on $N$ points.
Hence, when $N\geq n_\epsilon$, we have 
we trivially have $\mathrm{dist}(\dn, \qe)\to0$ because $\qe\in\mathcal D$. 
\wjsay{What does this mean? What is dist?}

As $\seqin{\alpha^*_i}$ attains the maximum of \eqref{eq:kmemp_prac},
\begin{multline}
\liminf_{N\to\infty}\sum_{i=1}^N\alpha_i^* l(z_i)\geq\sum_{i=1}^{K^{(n_\epsilon)}}{\gamma_i^{(n_\epsilon)}}{l({y_i^{(n_\epsilon)})}}\\
 = \int l(x) \ d\qe (x).
\label{eq:prop2}
\end{multline}
Combining \eqref{eq:prop1} and \eqref{eq:prop2}, we obtain the result
\wjsay{Why is this? I cannot follow this implication.}
$$
\lim_{N\to\infty}\sum_{i=1}^N\alpha_i^* l(z_i) = \int l(x) \ dP^*(x).
$$

\begin{remark}
    This proof is self-contained and is different from the the proof of \cite{kanagawaRecoveringDistributionsGaussiana} which requires $l$ to be a function in the Besov space.
\end{remark}
\end{todo}        
\end{proof}

\subsection{RKHS ambiguity set with known first $p$ moments}
The embedding associated with this kernel is
\begin{multline}
\hat\mu
= \int k(x , \cdot ) d P(x) 
= \int \left(x^\top (\cdot) + 1\right)^2\ d P(x)\\
= (\cdot) ^\top \mathbb{E}{x x^\top}  (\cdot) + 2 \mathbb{E} {x}^\top(\cdot) +1
\label{eq:pol_kme}
\end{multline}
By virtue of the kernel trick,
\begin{equation}
        \begin{aligned}
& \|\mu_P - \hat\mu\|_\rkhs = \|\sum_{i=1}^N\alpha_i\phi(x_i) - \int k(x , \cdot ) d P(x)\| \\
& = \alpha^\top K_z\alpha - 2 \sum_{i=1}^N\alpha_i\int k(x, z_i) dP(x) \\
&+ \int\int k(x, x') dP(x)dP(x')\\
& = \alpha^\top K_z\alpha - 2 \sum_{i=1}^N\alpha_i({z_i}^\top \mathbb{E}{x x^\top} z_i + 2 \mathbb{E} {x}^\top z_i +1) \\
&\quad + \mathrm{Tr}(\mathbb{E}{x x^\top}\mathbb{E}{x x^\top}  ) + 2 \mathbb{E} {x}^\top\mathbb{E} {x} +1.
        \end{aligned}
\end{equation}
This results in \eqref{eq:mp_pol}.

\section{Additional detailed experimental setup}\label{sec:additional}
\subsection{Setup for Sec.~\ref{sec:toy}}
\label{sec:app_exp1}
For the simulation to produce Fig.~\ref{fig:cantelli_n100}, we generate \(M=100\) data samples \(\{x_i\}_{i=1}^{M}\overset{i.i.d.}{\sim} N(0,1)\) (distribution unknown to the optimizer) and set $c=2.5$. To form the new expansion points $\seqin{z_i}$ in \eqref{eq:kme_indicator}, we first generate a set of $L=100$ grid points $\{y_i\}_{i=1}^L$ on $[0, 5]$, then we add the original data samples $\{x_i\}_{i=1}^{M}$ to the grid points to form the new expansion points. That is, the set of new expansion points $\seqin{z_i}$ consists of both the grid points and the original expansion points $\seqin{z_i}=\{y_i\}_{i=1}^L\cup\{x_i\}_{i=1}^{M}$, $N=200$ in total.

For the experiment in Fig.~\ref{fig:transport}, we add \(L=10\) new expansion points $\{y_i\}_{i=1}^L$ to the other side of \(c=2.5\) (by simply gridding the interval $[c,c+1]$). Together with the original expansion points \(\{x_i\}_{i=1}^M\), they form the set of expansion points $\seqin{z_i}$ for the new embedding $\kme{\alpha_i}{z_i}$.

\subsection{Setup for Sec.~\ref{sec:smpc}}
\label{sec:app_exp2}
For the scenario MPC simulation, we sample the i.i.d. uncertainty realizations \(\{s_1, \dots s_n\}\sim N(m,\Sigma)\), where \(m=[0.5\ \ 0]^{T},\Sigma=\left[\begin{array}{cc}  {0.01^{2}} & {0} \\  {0} & {0.1^{2}}  \end{array}\right]\). 
Then, the scenarios are propagated through the nonlinear dynamics.
The continuous-time dynamics is transcribed with numerical integration. 
The total time horizon is \(1.0\)s and we consider \(10\) control steps in this experiment. 

\end{document}